\documentclass[10pt]{amsart}
\usepackage{amssymb, amscd, amsmath, amsthm, epsfig, latexsym, enumerate}

\renewcommand{\geq}{\geqslant}
\renewcommand{\leq}{\leqslant}

\newcommand{\Z}{\mathbb Z}

\newtheorem{theorem}{Theorem}
\newtheorem{lemma}[theorem]{Lemma}

\newtheorem*{thm}{Theorem}
\newtheorem*{lem}{Lemma}
\newtheorem*{cor*}{Corollary}
\newtheorem*{add}{Addendum}

\begin{document}

\title{virtually solvable groups with balanced presentations}

\author{Jonathan A. Hillman }
\address{School of Mathematics and Statistics\\
     University of Sydney, NSW 2006\\
      Australia }

\email{jonathanhillman47@gmail.com}

\begin{abstract}
If a virtually solvable group with finite Hirsch length has a balanced presentation
then its abelianization has rank at most 3.
However the Hirsch lengths are essentially unrestricted.
We give some related examples of embeddings of 3-manifolds in $S^4$.
\end{abstract}


\maketitle
A group $G$ is \emph{homologically balanced} if $\beta_2(G;F)\leq\beta_1(G;F)$ 
for all prime fields $F$.
Groups with balanced presentations are homologically balanced.
In \S1 we show that if a group $G$ is homologically balanced,
virtually solvable and of finite Hirsch length then 
the abelianization $G^{ab}$ has rank at most 3.
(The argument extends one based on \cite{CEM} which was used for nilpotent groups in \cite{Hi22}.)
In \S2 we show that there are homologically balanced metabelian groups
with abelianization of rank 1, 2 or 3 and arbitrarily large Hirsch length.
We then consider 2-generator nilpotent groups, and show that no such group $G$ with 
$G/\gamma_4G\cong{F(2)/\gamma_4F(2)}$ can have a balanced 2-generator presentation.
(We use the notation of \cite{Rob} for the lower central series.)

Any embedding of a closed 3-manifold in $S^4$ has at least one
complementary region with homologically balanced fundamental group.
In \S4 we give examples in which the fundamental groups of each complementary region are solvable
(indeed, metabelian), but one is not homologically balanced.
(Details of the construction of such embeddings are confined to an appendix.)

\section{homologically balanced virtually solvable groups}

The crux of the argument for our main result (Theorem 1) is an appeal to work of
Charalambous, Evans and Miller on the Buchsbaum-Eisenbud Conjecture in commutative algebra.

\begin{lem}
\cite{CEM}
Let $R$ be a Cohen-Macaulay local ring of Krull dimension $d$ 
and with residue field $k$.
Let $M$ be a nonzero $R$-module of finite length.
Then
\[\dim_k\mathrm{Tor}_1^R(k,M)\geq\dim_kk\otimes_RM+d-1. \hfil\qed
\]
\end{lem}

In our applications $R$ shall be a localization of $k[x^{\pm1},y^{\pm1},z^{\pm1}]$,
for $k=\mathbb{Q}$ or $\mathbb{F}_p$.
We shall omit the coefficients $\mathbb{Q}$ from the notation
for homology and Betti numbers in Theorem 1,
as we use only rational homology there.

\begin{theorem}
Let $G$ be a finitely generated group with an epimorphism $f:G\to\Z^3$ such that
the abelianization $K^{ab}$ of $K=\mathrm{Ker}(f)$ has finite rank.
If $\beta_2(G)\leq\beta_1(G)$ then $\beta_2(G)=\beta_1(G)=3$.
\end{theorem}

\begin{proof}
Let $\Gamma=\mathbb{Q}[G/K]\cong\mathbb{Q}[x^{\pm1},y^{\pm1},z^{\pm1}]$,
and let $\varepsilon:\Gamma\to\mathbb{Q}$ be the augmentation homomorphism.
Conjugation in $G$ induces an action of $G/K$ on $K^{ab}$,
and gives $A=H_1(K)=\mathbb{Q}\otimes{K^{ab}}$ the structure of a $\Gamma$-module.
This has finite length since $K^{ab}$ has finite rank.
The LHS spectral sequence for $G$ as an extension of $G/K$ by $K$ gives rise to an exact sequence
\[
0\to{J}\to{H_2(G)}\to{H_2(G/K)}\to{H_0(G/K;A)}\to
{H_1(G)}\to{H_1(G/K)}\to0,
\]
where $J$ maps onto the cokernel of the differential $d^2_{3,0}:H_3(G/K)\to{H_1(G/K;A)}$.

Let $b_i=\dim_\mathbb{Q}H_i(G/K;A)$.
The exact sequence gives 
\[\dim_\mathbb{Q}J+3+\beta_1(G)=\beta_2(G)+b_0+3,
\]
since $H_2(G/K)\cong{H_1(G/K)}\cong\mathbb{Q}^3$.
Hence $\beta_2(G)=\beta_1(G)+\dim_\mathbb{Q}J-b_0$.
Moreover,  $\dim_\mathbb{Q}J\geq\max\{b_1-1,0\}$,
since $H_3(G/K)\cong\mathbb{Q}$.
Hence 
\[
\beta_2(G)\geq\beta_1(G)+\max\{b_1-1,0\}-b_0\geq\beta_1(G)-b_0.
\]

The terms $H_i(G/K;A)=Tor_i^\Gamma(\mathbb{Q},A)$ are finitely generated 
$\Gamma$-modules which are annihilated by $\mathfrak{m}=\mathrm{Ker}(\varepsilon)$,
since $\mathfrak{m}\mathbb{Q}=0$.
Let $R=\Gamma_\mathfrak{m}$ be the localization
of $\Gamma$ at the maximal ideal $\mathfrak{m}$.
The ring $R$ is a regular local ring of Krull dimension $d=3$,
which contains its quotient field $\mathbb{Q}$.
Since $\mathbb{Q}_\mathfrak{m}=\mathbb{Q}$ and localization is exact,
\[
H_i(G/K;A)\cong{Tor_i^\Gamma(\mathbb{Q},A)_\mathfrak{m}}=
Tor_i^R(\mathbb{Q},A_\mathfrak{m})\quad\mathrm{for~all}~i\geq0.
\]
Suppose that $H_0(G/K;A)\not=0$.
Since $b_0>0$ and regular local rings are Cohen-Macaulay,
$b_1\geq{b_0+d-1}=b_0+2$ \cite{CEM}.
Hence  $\beta_2(G)\geq\beta_1(G)+1$, which contradicts the hypothesis.
Therefore we must have $H_0(G/K;A)=0$,
and so ${H_1(G)}\cong{H_1(G/K)}$.
Hence $\beta_1(G)=3$ and $\beta_2(G)\geq\beta_1(G)$, 
and so $\beta_2(G)=\beta_1(G)=3$.
\end{proof}

\begin{cor*}
If $G$ is a virtually solvable group with Hirsch length $h(G)<\infty$
and such that $\beta_2(G)\leq\beta_1(G)$
then $\beta_1(G)\leq3$.
\end{cor*}

\begin{proof}
The theorem applies since the rank of any abelian section of $G$ is bounded by $h(G)$.
\end{proof}

The result of \cite{CEM} applies equally well when $k$ is a prime field $\mathbb{F}_p$.

\begin{add}
If $G$ has a balanced presentation and there is an epimorphism $f:G\to\mathbb{Z}^3$
with finitely generated kernel $K$ then $G^{ab}\cong\mathbb{Z}^3$ (and $K=G'$).
\end{add}

\begin{proof}
If $G$ has a balanced presentation and $K$ is finitely generated then
$\beta_2(G;\mathbb{F}_p)\leq\beta_1(G;\mathbb{F}_p)$ and 
$\mathbb{F}_p\otimes{K^{ab}}=K^{ab}/pK^{ab}$ is finite, 
for all primes $p$.
Hence if we repeat the argument with coefficients $\mathbb{F}_p$,
we find that $H_0(G/K;\mathbb{F}_p\otimes{K^{ab}})=0$ for all primes $p$.
It follows easily that $f$ induces an isomorphism on abelianizations.
\end{proof}

If $G$ is finitely generated and nilpotent then every subgroup is finitely generated.
Thus the addendum applies if $G$ is nilpotent, has a balanced presentation
and maps onto $\mathbb{Z}^3$.
In this case $\mathfrak{m}^kA=0$ for $k$ large,  
and so $A/\mathfrak{m}A=H_0(G/K;A)=0$ implies that $A=0$.
Similarly, $\mathbb{F}_p\otimes{K^{ab}}=0$ for all primes $p$.
Hence $K^{ab}=0$,  and so $K=1$, since $K$ is nilpotent.
Therefore either $G\cong\Z^3$ or $\beta_1(G)\leq2$ \cite{Hi22}.

In Theorem B of  \cite{Wi91} it is shown that if $G$ is a group with a presentation
on $n$ generators and $r$ relations and $G^{ab}$ can be generated by $d$ elements 
then either
\begin{enumerate}
\item $G$ has an infinite residually finite $p$-torsion quotient, for some prime $p$; or
\item$n-r\leq{d}-\frac14(d^2-1)$.
\end{enumerate}
In particular, if $G$ is virtually solvable and has a balanced presentation
then the second option must hold,
since finitely generated virtually solvable groups do not have infinite torsion groups as quotients.
Hence $d\leq4$ for such groups.
The above addendum has an additional finiteness condition,
but comparison of these results suggests the question:
can the bound on $d$ be improved at all in our context?

If a virtually solvable group has deficiency $>0$ then it is either $\Z$ or
a Baumslag-Solitar group $BS(1,m)$ with presentation 
$\langle{a,t}\mid{tat^{-1}=a^m}\rangle$ \cite[Corollary 2.6.1]{FMGK},\cite{Wi96}.

\section{homological balance and hirsch length}

In this section we show that there are homologically balanced metabelian groups realizing
every possible combination of Betti numbers and Hirsch length compatible with Theorem 1, 
other than $\beta_1=2$ and $\beta_2=0$.

We observe first that no group $G$ with $\beta_1(G)=2$ and $\beta_2(G)=0$ can be virtually solvable.
For if $\beta_1(G)=2$ then $G$ maps onto $\Z^2$, 
and if also $\beta_2(G)=0$ then any homomorphism from $F(2)$ to $G$ 
which lifts the abelianization $F(2)\to\Z^2$ must be injective,
since it induces isomorphisms on the quotients by the rational lower central series \cite{St65}.

The values of the pair $(\beta_1,\beta_2)$ with $\beta_2\leq\beta_1$
other than (2,0) allowed by Theorem 1 are realized by the free abelian groups $1,\Z,\Z^2,\Z^3$ 
and the semidirect product $\Z^2\rtimes_\Theta\Z$, where 
$\Theta=\left(\smallmatrix 0&-1\\
1&1\endsmallmatrix\right)$.
These groups realize the Hirsch lengths $\leq3$,
and all have balanced presentations.

Let $n\geq2$ and let $G_n$ be the group with the balanced presentation
\[
\langle {t,z}\mid{t^nzt^{-n}=tzt^{-1}z^{-1}}=z^{-1}tzt^{-1}\rangle.
\]
Let $z_i=t^izt^{-i}$, for $i\in\mathbb{Z}$.
Then $G_n^{ab}\cong\mathbb{Z}$ and $G_n'$ is generated by the conjugates $z_0,\dots,z_{n-1}$.
Conjugating the relation $z_1\leftrightharpoons{z_0}$ by powers of $t$ gives 
$z_{i+1}\leftrightharpoons{z_i}$, for all $i$.
Conjugating the relation $z_n=z_1z_0^{-1}$ gives $z_{n+1}=z_2z_1^{-1}$,
and so $z_2z_1^{-1}\leftrightharpoons{z_n}=z_1z_0^{-1}$.
Since $z_1\leftrightharpoons{z_0,z_2}$,
it follows that $z_2\leftrightharpoons{z_0}$.
Continuing in this way, we find that  $z_i\leftrightharpoons{z_j}$ for all
$0\leq{i,j}\leq{n}$, and so $G_n'$ is abelian.
It is then clear that $G_n'\cong\mathbb{Z}^n$ and so $G_n$ is metabelian and has Hirsch length $n+1$.
Thus the Hirsch length of a solvable group with a balanced presentation may be arbitrarily large.

Indeed, it may be infinite. 
These groups are closely related to the group $B\cong\mathbb{Z}^\infty\rtimes\mathbb{Z}^2$
with the 3-generator balanced presentation
\[
\langle{a,s,t}\mid{st=ts}, ~s^{-1}as=t^{-1}ata=at^{-1}at\rangle
\]
given by G. Baumslag in \cite{Ba72}.

If $G=K\rtimes\Z$ then the LHS spectral sequence for the homology of $G$ 
as an extension of $\Z$ by $K$ reduces to the Wang exact sequence
\begin{equation*}
\begin{CD}
\to{H_2(K;\Z)}@>t-1>>{H_2(K;\Z)}\to{H_2(G;\Z)}\to{H_1(K;\Z)}@>t-1>>{H_1(K;\Z)}\to
\end{CD}
\end{equation*}
\[
\\ \to{H_1(G;\Z)}\to\Z\to0.
\]
If $G^{ab}\cong\Z$ then $t-1$ acts invertibly on $H_1(K;\Z)$,
and so $H_2(G;\Z)$ is the quotient of $H_2(K;\Z)$ by the action of $t-1$.
If, moreover, $K$ is abelian then $H_2(K;\Z)\cong{K\wedge{K}}$,
with the diagonal action of $G/K$.
When $G=G_2$ and $K=G_2'\cong\Z^2$ this action is trivial, and so $H_2(G_2;\Z)\cong\Z$.
However, when $n>2$ it turns out that $t-1$ acts invertibly on $H_2(G_n';\Z)$
and so $H_2(G_n;\Z)=0$.
(This discrepancy becomes clearer if we extend the coefficients to $\mathbb{C}$
and diagonalize $t$, since the characteristic polynomial for the action of $t$ on $G_n'$
is $t^n-t+1$, which has a reciprocal pair of roots only when $n=2$.)

It follows from the K\"unneth theorem that $\Z\times{G_n}$ 
is homologically balanced (with $\beta_1=2$), for all $n\geq2$,
while $\Z^2\times{G_n}$ is homologically balanced (with $\beta_1=3$), for all $n>2$.
The remaining possibilities with $\beta_1=\beta_2=3$ and $h=4$ or 5 are realized by $\Z^2\times{H}$, 
where $H$ is $BS(1,2)$ or the group with presentation
$\langle{a,b,t}\mid{ab=ba},~tat^{-1}=a^2,~tbt^{-1}=b^2\rangle$,
respectively.
In the appendix we shall show that there are no homologically balanced poly-$\Z$ groups 
with $\beta_1=3$ and $h=4$ or 5.

The groups with the balanced presentations $\langle{x,y}\mid [x,[x,[x,y]]]=[y,[x,y]]=1\rangle$ and 
$\langle {t,z}\mid{t^4zt^{-2}z^{-1}t^{-1}zt^{-1}z^{-1}=[z,tzt^{-1}]=1}\rangle$
are metabelian and have $\beta_1=2$.
Their Hirsch lengths are 4 and 5, respectively.
We do not know whether there are any solvable groups with balanced presentations,
$\beta_1=2$ or 3,
and Hirsch length greater than 5.

\section{2-generator nilpotent groups}

Let $G$ be a finitely generated group.
For each $r\geq2$, 
let $\Phi_r(G)$ be the kernel of the homomorphism 
from $H_2(G)$ to $H_2(G/\gamma_{r-1}G)$ induced by the projection. 
Thus $\Phi_2(G)=H_2(G)$, and $\Phi_{r+1}(G)\leq\Phi_r(G)$, for all $r$.
These subgroups define the \emph{Dwyer filtration} of $H_2(G;\mathbb{Q})$.
If $f:\sigma\to\pi$ is a group homomorphism such that $H_1(f)$ is an isomorphism 
then $f$ induces an isomorphism $H_2(\sigma)/\Phi_r(\sigma)\cong{H_2(\pi)/\Phi_r(\pi)}$
if and only if it induces an isomorphism 
$\sigma/\gamma_r\sigma\cong\pi/\gamma_r\pi$ \cite[Theorem 1.1]{Dw75}.
The proof is by induction on $r$, based on an exact sequence
\[
0\to\Phi_k(G)/\Phi_{k+1}(G)\to\Phi_k(G/\gamma_kG)\to
\mathbb{Q}\otimes\gamma_kG/\gamma_{k+1}G\to0
\]
which follows from a diagram chase \cite[Lemma 1.3]{Dw75}.

Let $F(\beta)$ be the free group of rank $\beta$.
Then $\Phi_r(G)=H_2(G)$ if and only if $G/\gamma_rG\cong{F(\beta)/\gamma_rF(\beta)}$. 
If $G=F(\beta)/\gamma_kF(\beta)$ then
$H_2(G;\mathbb{Z})\cong\gamma_kF(\beta)/\gamma_{k+1}F(\beta)$,
by the five-term exact sequence of low degree for the homology of
$G$ as a quotient of $F(\beta)$.
This abelian group has rank $r_k(\beta)=\frac1k\Sigma_{d|k}\mu(d)\beta^{\frac{k}d}$,
where $\mu$ is the M\"obius function, 
by the Witt formulae \cite[Theorems 5.11 and 5.12]{MKS}.
Hence the only relatively free nilpotent groups $G$ with $\beta_2(G)\leq\beta_1(G)$ are  
$G\cong\mathbb{Z}^\beta$ with $\beta\leq3$ and $G\cong{F(2)/\gamma_3F(2)}$.
These groups have balanced presentations and Hirsch length $\leq3$.
We also see that $h(F(2)/\gamma_4F(2))=5$ and $h(F(2)/\gamma_5F(2))=8$.

The following assertion summarizes the main results of work of Freedman,  Hain and Teichner,
based on the Dwyer filtration and an analogue of Roquette's proof of the Golod-Shafarevich inequality
 \cite{FHT}.
Let  $p_k$ be the rank of $\Phi_k(G)/\Phi_{k+1}(G)$.
If $G$ is nilpotent of class $n$ then $p_r=0$ for $r>n+1$,
and $\Sigma_{r\geq2}p_r=\beta_2(G)$.

\begin{thm}
 \cite[Theorem 2 and Proposition 5]{FHT}
If $G$ is a finitely generated nilpotent group such that $h(G)>2$ then 
$p(t)=1-\beta_1(G){t}+\Sigma_{r\geq2}p_rt^r>0$ for $0<t<1$.

Hence if $\Phi_r(G)=H_2(G)$ then 
$\beta_2(G)t^r-\beta_1(G){t}+1>0$ for $0<t<1$, 
and so $\beta_2(G)>\frac{(r-1)^{r-1}}{r^r}\beta_1(G)^r.$
\qed
\end{thm}

Setting $r=2$ in this theorem gives $\beta_2(G)>\frac14\beta_1(G)^2$,
since $\Phi_2(G)=H_2(G)$ for any group $G$.
Hence if $G$ is nilpotent and $\beta_2(G)=\beta_1(G)$ then $\beta_1(G)\leq3$.
In \cite{Hi22} we showed that $\mathbb{Z}^3$ is the only homologically balanced nilpotent group
with $\beta_1(G;\mathbb{Q})=3$.
Thus the remaining interest is in nilpotent quotients of $F(2)$.
If $G/\gamma_5G\cong{F(2)/\gamma_5F(2)}$ 
then setting $r=5$ gives $\beta_2(G)>\beta_1(G)=2$.
We may extract a little more from \cite[Proposition 5]{FHT}
when $G/\gamma_4G\cong{F(2)/\gamma_4F(2)}$.

\begin{lemma}
If $G$ is a nilpotent group such that $G/\gamma_4G\cong{F(2)/\gamma_4F(2)}$ and $\beta_2(G)=\beta_1(G)$ then $h(\gamma_4G/\gamma_5G)=1$.
\end{lemma}

\begin{proof}
If $G$ is such a group  then $p_2=p_3=0$ and $\Sigma_{r\geq2}p_r=2$.
If $p_4\leq1$ then $p(t)\leq{t^5+t^4-2t+1}$ for $0<t<1$,
and so  $p(0.65)<0$.
This contradicts \cite[Proposition 5]{FHT},
and so $p_4=2$.
Since $\Phi_4(G/\gamma_4G)=H_2(G/\gamma_4G)$ has rank 3,  
it follows from \cite[Lemma 1.3]{Dw75} that 
$\gamma_4G/\gamma_5G$ has rank $1$.
\end{proof}

The authors of \cite{FHT} conjecture that there are no such nilpotent
groups $G$ (with $G/\gamma_4G\cong{F(2)/\gamma_4F(2)}$ and
$\beta_2(G)=\beta_1(G)$).
We shall show that  no such nilpotent group can have a balanced presentation
with 2 generators and 2 relations.

\begin{lemma}
If $\beta\geq2$ and $k>2$ then $r_k(\beta)=\frac1k\sum_{s|k}\mu(s)\beta^{\frac{k}s}>2k$.
\end{lemma}

\begin{proof}
The assertion is easily checked by hand for $k\leq9$,
and so we may assume that $k\geq10$.
Let $\omega(k)$ be the number of distinct prime divisors of $k$. 
Then $\omega(k)<k$ for all $k$.
The sum $\sum_{s|k}\mu(s)\beta^{\frac{k}s}$
may be expanded as the  sum of $\omega(k)+1$ terms
$\beta^k+\sum_{1<i\leq\omega(k)}(-1)^i\sigma_i$,
where $\sigma_i$ is the sum of all terms $\beta^{\frac{k}s}$
with $s$ a product of $i$ distinct prime divisors of $k$.
Now $\sigma_1=\sum_{p|k}\beta^{\frac{k}p}\leq\omega(k)\beta^{\frac{k}2}$,
and so $\sigma_1<\beta^k$.
Similarly, 
$\sigma_{i+1}=\sum_{p|k}\sum_{i,\not{p}}(\beta^{\frac1p})^{\frac{k}s}$,
where the inner sum is taken over products of $i$ prime factors of $k$ which are 
distinct and \emph{not} $p$.
A simple induction shows that $\sigma_{i+1}<\sigma_i$ for all $1<i<\omega(k)$.
Hence $\sigma_{2j}-\sigma_{2j+1}>0$ if $2j<\omega(k)$, and so
\[
kr_k(\beta)>\beta^k-\sigma_1\geq\beta^k-\omega(k)\beta^{\frac{k}2}>
\beta^k-k\beta^{\frac{k}2}=\beta^{\frac{k}2}(\beta^{\frac{k}2}-k).
\]
We have $\beta^\frac{k}2\geq2k$, since $k\geq10$,
and so $r_k(\beta)>\frac1k2k(2k-k)>2k$.
\end{proof}

We shall only need this estimate for $\beta=2$.

\begin{theorem}
If $G$ is a group such that $G/\gamma_4G\cong{F(2)/\gamma_4F(2)}$ and 
which has a balanced presentation with $2$ generators and $2$ relations
then $G$ is not nilpotent.
\end{theorem}

\begin{proof}
We may suppose that $G=F/R$, where $F=F(x,y)$ and $R$ is normally generated by
$u,v\in\gamma_4F$.
The image of $u$ and $v$ in $\gamma_4(F)/\gamma_5(F)$ generates a subgroup
of rank at most 2, and hence the subgroup of $\gamma_5(F)/\gamma_6(F)$
generated by the commutators $[x,u],[x,v],[y,u]$ and $[y,v]$ has rank at most 4.
Continuing in this way, we find that the rank of the image of higher commutators of $u$
and $v$ in  $\gamma_k(F)/\gamma_{k+1}(F)$ is at most $2(k-3)$.
Hence there is no value of $k$ such that $\gamma_k(F)\leq{R}$, and so $G$ is not nilpotent.
\end{proof}

If $G$ is a nilpotent group which can be generated by $d$ elements and has 
a balanced presentation then does it have a presentation with $d$ generators and $d$ relations?

The group with presentation $\langle{x,y}\mid [x,[x,[x,y]]]=[y,[x,y]]=1\rangle$
is the only torsion-free nilpotent group of Hirsch length 4 with a balanced presentation,
and there are none such of Hirsch length 5 \cite{Hi22}.
Are there good upper bounds for the deficiency of a virtually nilpotent group in terms of its Hirsch length?

\section{an application to embedding 3-manifolds in $S^4$}

Let  $M$ be a closed 3-manifold with an embedding $j:M\to{S^4}$ 
with complementary regions $X$ and $Y$.
Then $X$ and $Y$ are connected, $\chi(X)+\chi(Y)=\chi(S^4)=2$, 
$\beta_{3-i}(Y)=\beta_i(X)$ and
$\beta_i(X)+\beta_i(Y)=\beta_i(M)$ for $i=1,2$ and $\beta_j(X)=\beta_j(Y)=0$ for $j>2$.
We may assume that $\chi(X)\leq1\leq\chi(Y)$,
and then $\pi_X$ is homologically balanced. 
(These assertions summarise \cite[Lemma 2.2]{Hi}.)

In  \cite[Theorem 9.7]{Hi}  we show that if $\pi_X=\pi_1(X)$ and $\pi_Y=\pi_1(Y)$
are virtually solvable then $\beta_1(M)\leq8$. 
The full statement of this theorem follows:

\begin{thm}
\cite[Theorem 9.7]{Hi}
Let $j:M\to{S^4}$ be a restrained embedding determined by a bipartedly trivial link.
If $\beta=\beta_1(M;\mathbb{Q})$ is odd then $\chi(X)=0$ and $\chi(Y)=2$,
and $X$ is aspherical.
If, moreover,  $\pi_X$ is almost coherent or elementary amenable then
$\pi_X\cong\mathbb{Z}$ or $BS(1,m)$, for some $m\not=0$,
and $\beta=1$ or $3$.

If $\beta$ is even then $\chi(X)=\chi(Y)=1$ and $\pi_X$ and $\pi_Y$ have balanced presentations.
Moreover, if $\pi_X$ and $\pi_Y$ are virtually solvable then $\beta=0$, $2$, $4$, $6$ or $8$.
\qed
\end{thm}  

Here ``$j$ is restrained" means ``$\pi_X$ and $\pi_Y$ have no non-cyclic free subgroups".
The assumption that $j$ is determined by a bipartedly trivial link is used only to ensure that
$X$ and $Y$ are homotopy equivalent to finite 2-complexes,
while the bound 8 derives from the above-mentioned consequence of \cite[Theorem B]{Wi91}.

Suppose that $\pi_X$ and $\pi_Y$ are virtually solvable and of finite Hirsch length.
If $\beta_1(X)=\beta_1(\pi_X)>3$ then $\beta_2(X)\geq\beta_2(\pi_X)>\beta_1(X)$, 
by Theorem 1, and so $\chi(X)\geq2$.
Similarly,  $\chi(Y)\geq2$, since we then also have $\beta_1(Y)=\beta_2(X)>3$.
Hence $\chi(X)+\chi(Y)\geq4$, which is wrong.
Therefore $\beta_1(X)\leq3$ and $\beta_1(Y)\leq3$.
Moreover,  if $\beta_1(X)<3$ then $\beta_1(Y)<3$ also, and so
$\beta_1(M)\leq4$ or $\beta_1(M)=6$.

This sharpens \cite[Theorem 9.7]{Hi} in two respects; 
the embedding need not be determined by a bipartedly trivial link,
and there are simple examples (with $\pi_X$ and $\pi_Y$ abelian!) realizing the remaining
allowed values for $\beta_1(M)$.
On the other hand, we need finiteness of the Hirsch length for the appeal to Theorem 1.

In general, if $\chi(X)\leq1$ then $\pi_X$ is homologically balanced \cite[Lemma 2.2]{Hi}.
However,  it  is not generally true that both complementary regions 
should have fundamental groups with homologically balanced presentations.
If $\pi_X$ is virtually solvable and $\chi(X)<1$ then $\beta_1(\pi_X)=\beta_1(X)>\beta_2(X)$,
and $\beta_2(X)\geq\beta_2(\pi_X)$, by Hopf's Theorem.
Hence $\beta_1(\pi_X)=\beta_2(\pi_X)+1$ (by the observations deriving from \cite{St65}
following the corollary above),
and so $\chi(X)=0$ and $\beta_1(\pi_X)=1$ or 2.
Hence  $\beta_2(Y)=\beta_1(Y)+1$ and $\beta_2(Y)=1$ or 0.

In order to show that there are indeed such examples,
we shall appeal to a result of Quinn which asserts that a pair of 2-complexes $C$ and $D$
which satisfy the conditions $H_1(D;\mathbb{Z})\cong{H^2(C;\mathbb{Z})}$
and $H_2(D;\mathbb{Z})\cong{H^1(C;\mathbb{Z})}$ derived from Alexander duality in $S^4$
can be realized up to homotopy as the complementary regions of 
an embedding of a 3-manifold in $S^4$ \cite{Qu01}.

Let $C$ and $D$ be the 2-complexes corresponding to the presentations
\[
\langle{a,s}\mid{sas^{-1}=a^7}\rangle
\]
and
\[
\langle{b,v}\mid{b^6=v^3=(bv)^2=1}\rangle,
\]
of the solvable Baumslag-Solitar group $BS(1,7)$ and the 2-dimensional crystallographic group 
$[p6]=\pi^{orb}S^2(2,3,6)\cong\mathbb{Z}^2\rtimes{C_6}$, 
respectively.
The  complex $C$ is aspherical,
$H_1(C;\mathbb{Z})\cong\mathbb{Z}\oplus\mathbb{Z}/6\mathbb{Z}$,
and $H_2(C;\mathbb{Z})=0$,
while $H_1(D;\mathbb{Z})\cong\mathbb{Z}/6\mathbb{Z}$ and $H_2(D;\mathbb{Z})\cong\mathbb{Z}$.
It follows easily from the Universal Coefficient Theorem that the pair $C,D$ satisfies the criterion of
\cite{Qu01}.
The group  $\pi_1D$ is not homologically balanced,
since $H_1(\pi_1D;\mathbb{Q})=0$ and
$H_2(\pi_1D;\mathbb{Q})\cong{H_0(C_6;H_2(\Z^2;\mathbb{Q}))}\cong\mathbb{Q}$.
(Is there a similar example with $X\simeq{S^1}$, 
$\pi_Y$ virtually solvable and $\beta_2(\pi_Y;\Z)=1$?)

In our second example $\pi_X\cong\Z^2$ and $\pi_Y$ is a semidirect product
$\Z[\zeta]\rtimes\Z$, where $\zeta$ is a primitive 10th root of unity.
The minimal polynomial for $\zeta$ is
\[
\phi_{10}(t)=t^4-t^3+t^2-t+1=(t-\zeta)(t-\zeta^3)(t-\zeta^7)(t-\zeta^9).
\]
Let $G_\zeta$ be the group with presentation
\[
\langle{t,z}\mid{t^4zt^{-1}z^{-1}t^{-1}zt^{-1}z^{-1}t^{-1}z}=[t^2zt^{-2},z]=[tzt^{-1},z]=1\rangle.
\]
Then $G_\zeta^{ab}\cong\Z$,  so $\beta_1(G_\zeta)=1$, and $G_\zeta'\cong\Z[\zeta]$.
The Wang sequence of homology for $G_\zeta$ as an extension of $\Z$ by $G_\zeta'$
gives an exact sequence
\begin{equation*}
\begin{CD}
{H_2(G_\zeta';\Z)}@>t-1>>{H_2(G_\zeta';\Z)}\to{H_2(G_\zeta;\Z)}\to
{H_1(G_\zeta';\Z)}@>t-1>>{H_1(G_\zeta';\Z)}\to0,
\end{CD}
\end{equation*}
since $H_1(G_\zeta;\Z)\cong\Z$.
Hence $t-1$ maps $H_1(G_\zeta';\Z)=G_\zeta'$ onto itself, and so it is an automorphism of $G_\zeta'$.
Since $G_\zeta'$ is abelian, $H_2(G_\zeta';\Z)\cong{G_\zeta'\wedge{G_\zeta'}}\cong\Z^6$, 
and the eigenvalues of the natural action of $t$ are the six
pairwise products of the distinct primitive 10th roots of unity.
Two of these products are 1.
Hence $H_2(G_\zeta';\Z)/(t-1)H_2(G_\zeta';\Z)$ has rank 2,
and so $\beta_2(G_\zeta)=2$.

Let $C=T$ and let $D$ be the finite 2-complex corresponding to the above presentation
of $G_\zeta$.
Then the pair $C,D$ satisfies the criterion of \cite{Qu01},
and so there is an embedding of a 3-manifold in $S^4$ with complementary regions 
$X\simeq{C=T}$ and $Y\simeq{D}$.
(In the appendix we shall outline a simple direct construction of such an embedding, 
avoiding appeal to the general argument of \cite{Qu01}.)

\section*{Appendix}

1.   We shall outline why a poly-$\Z$ group $G$ with $\beta_1(G;\mathbb{F}_2)\geq3$
and $h=4$ or 5 cannot be homologically balanced.
If $h(G)=4$  then $G$ is a Poincar\'e duality group of dimension 4 and $\chi(G)=0$,
so $\beta_2(G;\mathbb{F}_2)=2(\beta_1(G;\mathbb{F}_2)-1)>\beta_1(G;\mathbb{F}_2)$.

If $h(G)=5$ then $G$ has an abelian normal subgroup $A\cong\Z^2$, 
and the image of $G$ in $Aut(A)\cong{GL(2,\Z)}$ is either $\Z$,
$\Z\oplus\langle-I\rangle$, cyclic of order $1, 2,3,4$ or 6, 
or $C_2^2$.
Hence $G\cong{C}\rtimes\Z$, 
where $C$ is either nilpotent and has centre of rank $\geq2$
or $C\cong\Z\times(\Z^2\rtimes_{-I}\Z)$.  
Case-by-case Wang sequence arguments then show that $\beta_2(G)\geq4$.

2. The appeal to \cite{Qu01} may be avoided by giving a suitable 0-framed
bipartedly trivial link,  as in many of the examples in \cite{Hi}.
The meridians of each sublink represent the generators of one of the presentations,
and the corresponding component represents a relator for the other presentation.
We first add redundant generators and relators, if necessary,
so that each generator of the first presentation may be paired with a relator for the second presentation, and vice versa.
The exponent sum for each generator in a given relator
is the linking number of the corresponding components of the first and second sublinks,
and we may need to change the relators, 
to ensure symmetry of the linking numbers.

For instance, when the groups are $BS(1,7)$ and $[p6]=\pi^{orb}S^2(2,3,6)$ 
the simplest such pair of presentations that we have found are
\[
\langle{a,s,t}\mid{sa^{-1}s^{-1}a^7t^2=t=1}\rangle
\]
with 3 generators and 2 relators $B=sa^{-1}s^{-1}a^7t^2$ and $V=t$, and
\[
\langle{b,v}\mid{b^6=b^{-6}((bv)^2v^{-3})^3v^3=v^3(bv)^{-2}=1}\rangle.
\]
with 2 generators and 3 relators $A=b^6$,
$S=b^{-6}((bv)^2v^{-3})^3v^3$ and $T=v^3(bv)^{-2}$.
In this case $X\simeq{K(BS(1,7),1)}$ is aspherical.
However, we shall not give further details,
as the resulting link diagram is too messy to be easily drawn or enlightening.

There is a related example (with $X\simeq{S^1}$) based on the presentations
\[
\langle{u,v,w}\mid(uvw)^2vw^{-2}=(uvw)^2(w^2v^{-1})^2w=1\rangle
\]
for $\Z$ and 
\[
\langle{x,y}\mid{x^2=y^3=(xy)^7=1}\rangle
\]
for the hyperbolic triangle group $\pi^{orb}(S^2(2,3,7)$.
However the latter group is not virtually solvable,
and the associated link diagram is again too messy.
Is there a perfect virtually solvable group that we could use instead of 
$\pi^{orb}(S^2(2,3,7)$?

The pair of presentations for $\Z^2$ and $G_\zeta$ given in \S4 above leads to 
a neater example.
Here the only modification needed is the addition of a redundant generator and relator for $\Z^2$,
so that it has the presentation
\[
\langle{u,v,w}\mid{[u,v]=w=1}\rangle
\]
with 3 generators and 2 relators.
Then the only relators in this pair of presentations with 
a nonzero exponent sum are $Z=w$ and $W=t^4zt^{-1}z^{-1}t^{-1}zt^{-1}z^{-1}t^{-1}z$.

To realize this pair we construct a link $L$ with 5 components $U,V,W,T,Z$ and
such that the sublinks $\{T,Z\}$ and $\{U,V,W\}$ are trivial,
$\mathrm{lk}(W,Z)=1$ and all other linking numbers are 0.
Let $t,u,v,w,z$ be meridians for these components.
We choose $U, V,W$ to represent 
$t^4zt^{-1}z^{-1}t^{-1}zt^{-1}z^{-1}t^{-1}z$, $[t^2zt^{-2},z]$ and $[tzt^{-1},z]$
(respectively) in the free group $\pi_1(S^3\setminus\{T,Z\})$.
We may need to modify the position of the sublink $\{U,V,W\}$ so that $t$ and $z$ represent $[v,w]$ 
and $u$ (respectively) in the free group $\pi_1(S^3\setminus\{U,V,W\})$.
Let $M(L)$ be the result of 0-framed surgery on $L$ in $S^3$.
If we embed $S^3$ as the equator of $S^4$ then we may add 2-handles 
along $\{T,Z\}$to one hemisphere and along $\{U,V,W\}$ to the other,
to give an embedding of the 3-manifold $M(L)$ into $S^4$ with the desired properties
\cite[Chapter 2]{Hi}.

\setlength{\unitlength}{1mm}
\begin{picture}(100,86)(0.7,4)

\put(15,45){$T$}
\put(103,45){$Z$}
\put(38,65){$U$}
\put(38,23){$V$}
\put(38,40){$W$}

\put(33,14.2){$\triangleright$}
\put(31,12){$t$}
\put(88,14.2){$\triangleright$}
\put(86,12){$z$}
\put(59,20.1){$\gets$}
\put(51,30.1){$\gets$}
\put(59,47.1){$\gets$}

\linethickness{1pt}

\put(20,20){\line(0,1){60}}
\put(45,21.5){\line(0,1){2}}
\put(45,24.5){\line(0,1){4}}
\put(45,29.7){\line(0,1){1}}
\put(45,31.5){\line(0,1){4.5}}
\put(45,37){\line(0,1){7.5}}

\put(45,46){\line(0,1){1.5}}
\put(45,48.5){\line(0,1){1.7}}
\put(45,51.4){\line(0,1){1.6}}
\put(45,54.4){\line(0,1){1.6}}
\put(45,57){\line(0,1){7}}
\put(45,65.5){\line(0,1){3}}
\put(45,70.3){\line(0,1){3}}
\put(45,75.1){\line(0,1){3}}

\put(75,21.5){\line(0,1){3}}
\put(75,25.5){\line(0,1){12}}
\put(75,40.5){\line(0,1){16.5}}
\put(75,58.7){\line(0,1){7.5}}
\put(75,67.5){\line(0,1){1.4}}
\put(75,70.3){\line(0,1){5.5}}

\put(100,20){\line(0,1){60}}
\put(25,15){\line(1,0){15}}
\put(25,85){\line(1,0){15}}
\put(80,15){\line(1,0){15}}
\put(80,85){\line(1,0){15}}
\qbezier(20,20)(20,15)(25,15)
\qbezier(40,15)(45,15)(45,20)
\qbezier(75,20)(75,15)(80,15)
\qbezier(95,15)(100,15)(100,20)
\qbezier(20,80)(20,85)(25,85)
\qbezier(40,85)(45,85)(45,80)
\qbezier(75,80)(75,85)(80,85)
\qbezier(95,85)(100,85)(100,80)

\thinlines
\put(46,57.8){\line(1,0){30}}
\put(44,48){\line(1,0){30}}
\put(76,48){\line(1,0){3}}
\qbezier(79,48)(84,48)(84,53)
\qbezier(79,79)(84,79)(84,74)
\put(84,53){\line(0,1){21}}

\qbezier(43.3,48.7)(43.4,48)(44,48)
\qbezier(43.3,48.7)(43.4,49.4)(44,49.4)
\qbezier(45.5,49.4)(46.2, 49.4)(46.2,50.1)
\qbezier(45.5,50.8)(46.2, 50.8)(46.2,50.1)
\qbezier(43.3,51.5)(43.3,50.8)(44,50.8)
\qbezier(43.3,51.5)(43.3,52.2)(44,52.2)
\put(44,50.8){\line(1,0){1.5}}
\qbezier(45.5,52.2)(46.2, 52.2)(46.2,52.9)
\qbezier(45.5,53.6)(46.2, 53.6)(46.2,52.9)
\put(44,53.6){\line(1,0){1.5}}
\qbezier(43.3,54.3)(43.3,55)(44,55)
\qbezier(43.3,54.3)(43.3,53.6)(44,53.6)
\qbezier(45.5,55)(46.2, 55)(46.2,55.7)
\qbezier(45.5,56.4)(46.2, 56.4)(46.2,55.7)
\put(44,56.4){\line(1,0){1.5}}
\qbezier(43.3,57.1)(43.3,57.8)(44,57.8)
\qbezier(43.3,57.1)(43.3,56.4)(44,56.4)

\put(44,79){\line(1,0){35}}
\put(46,76.6){\line(1,0){30}}
\put(44,74.2){\line(1,0){30}}
\put(46,71.8){\line(1,0){28}}
\put(44,69.4){\line(1,0){32}}
\put(46,67){\line(1,0){30}}
\put(44,64.6){\line(1,0){30}}
\put(46,62.2){\line(1,0){28}}

\qbezier(42.8,63.4)(42.8,62.2)(44,62.2)
\qbezier(42.8,63.4)(42.8,64.6)(44,64.6)
\qbezier(42.8,68.2)(42.8,67)(44,67)
\qbezier(42.8,68.2)(42.8,69.4)(44,69.4)

\qbezier(42.8,73)(42.8,71.8)(44,71.82)
\qbezier(42.8,73)(42.8,74.2)(44,74.2)

\qbezier(42.8,77.8)(42.8,76.6)(44,76.6)
\qbezier(42.8,77.8)(42.8,79)(44,79)

\qbezier(76,57.8)(78,57.8)(78.2,60)
\qbezier(76,62.2)(78,62.2)(78.2,60)
\qbezier(76,64.6)(77.2,64.6)(77.2,65.8)
\qbezier(76,67)(77.2,67)(77.2,65.8)
\qbezier(76,69.4)(77.2,69.4)(77.2,70.6)
\qbezier(76,71.8)(77.2,71.8)(77.2,70.6)

\qbezier(76,74.2)(77.2,74.2)(77.2,75.4)
\qbezier(76,76.6)(77.2,76.6)(77.2,75.4)

\put(44,45){\line(1,0){30}}
\put(46,43){\line(1,0){28}}
\put(59.6,40){\line(1,0){16.4}}
\put(59.2,38){\line(1,0){16.8}}

\qbezier(43,44)(43,45)(44,45)
\qbezier(43,44)(43,43)(44,43)

\qbezier(76,43)(77.5,43)(77.5,41.5)
\qbezier(76,40)(77.5,40)(77.5,41.5)
\qbezier(76,45)(79.5,45)(79.5,41.5)
\qbezier(76,38)(79.5,38)(79.5,41.5)

\put(44,36.4){\line(1,0){30}}
\qbezier(43,35.6)(43,36.4)(44,36.4)
\qbezier(43,35.6)(43,34.8)(44,34.8)

\put(46,27.2){\line(1,0){28}}
\put(46,25){\line(1,0){30}}
\put(44,21){\line(1,0){32}}

\qbezier(43.4,21.6)(43.4,21)(44,21)
\qbezier(43.4,21.6)(43.4,22.2)(44,22.2)
\qbezier(43.4,24.4)(43.4,25)(44,25)
\qbezier(43.4,24.4)(43.4,23.8)(44,23.8)
\put(44,23.8){\line(1,0){1.5}}
\qbezier(45.5,23.8)(46.3, 23.8)(46.3,23)
\qbezier(45.5,22.2)(46.3, 22.2)(46.3,23)

\qbezier(43,28.2)(43,27.2)(44,27.2)
\qbezier(43,28.2)(43,29.2)(44,29.2)
\put(44,29.2){\line(1,0){22}}
\put(46,34.8){\line(1,0){20}}

\qbezier(66,29.2)(68,29.2)(68,31.2)
\qbezier(66,34.8)(68,34.8)(68,32.8)
\put(68,31.2){\line(0,1){1.6}}
\put(44,31){\line(1,0){10.8}}
\put(46,32.6){\line(1,0){9}}
\qbezier(55,31)(58.2,31)(58.2,34.2)
\qbezier(55,32.6)(56.6,32.6)(56.6,34.2)

\qbezier(56.6,37)(56.6,40)(59.6,40)
\qbezier(58.2,37)(58.2,38)(59.2,38)

\qbezier(43,31.8)(43,31)(43.8,31)
\qbezier(43,31.8)(43,32.6)(43.8,32.6)

\qbezier(76,25)(77.1,25)(77.1,26.1)
\qbezier(76,27.2)(77.1,27.2)(77.1,26.1)
\qbezier(76,21)(80,21)(80,25)
\qbezier(76,36.4)(80,36.4)(80,32.4)
\put(80,25){\line(0,1){7.4}}

\put(38,8){The bipartedly trivial link $L$.}

\end{picture}

{\bf Tool use.}
No AI tools were used in any part of this work.

%\newpage

\end{document}